RESEARCH ARTICLE

# Predicting the Last Zero of Brownian Motion with Drift

J. du Toit, G. Peskir & A. N. Shiryaev*

Given a standard Brownian motion $B^\mu = (B^\mu_t)_{0 \le t \le T}$ with drift $\mu \in I\!R$ and letting $g$ denote the last zero of $B^\mu$ before $T$, we consider the optimal prediction problem

$$V_* = \inf_{0 \le \tau \le T} \mathsf{E}|g - \tau|$$

where the infimum is taken over all stopping times $\tau$ of $B^\mu$. Reducing the optimal prediction problem to a parabolic free-boundary problem and making use of local time-space calculus techniques, we show that the following stopping time is optimal:

$$\tau_* = \inf \{ t \in [0,T] \mid B^\mu_t \le b_-(t) \text{ or } B^\mu_t \ge b_+(t) \}$$

where the function $t \mapsto b_-(t)$ is continuous and increasing on $[0,T]$ with $b_-(T) = 0$, the function $t \mapsto b_+(t)$ is continuous and decreasing on $[0,T]$ with $b_+(T) = 0$, and the pair $b_-$ and $b_+$ can be characterised as the unique solution to a coupled system of nonlinear Volterra integral equations. This also yields an explicit formula for $V_*$ in terms of $b_-$ and $b_+$. If $\mu = 0$ then $b_- = -b_+$ and there is a closed form expression for $b_\pm$ as shown in [10] using the method of time change from [4]. The latter method cannot be extended to the case when $\mu \neq 0$ and the present paper settles the remaining cases using a different approach.

*Supported by MIMS as Distinguished Visitor during June-July 2007.

## 1. Introduction

Imagine a stochastic process $X = (X_t)_{0 \le t \le T}$ observed over the time interval $[0,T]$, and let $g$ denote the last zero of $X$ before $T$. Formally, we set $g = \sup\{ t \in [0,T] \mid X_t = 0 \}$ so that $X_g = 0$ when $X$ is left-continuous or $X_{g-} = 0$ when $X$ has left limits (and the set is non-empty in either of the two cases). Clearly, at any time $t < T$ the value of $g$ is unknown (unless $X$ is trivial), and it is only at the terminal time $T$ itself that we know when the last zero $g$ of $X$ occurred. However, this is often too late: typically one wants to know how close $X$ is to $g$ at any time $t < T$ and then take some action based on this information. In the present paper we take up this question when $X$ is a standard Brownian motion $B^\mu = (B^\mu_t)_{0 \le t \le 1}$ with drift $\mu \in I\!R$ and search for a stopping time $\tau_*$ of $B^\mu$ that is as 'close' as possible to $g$. The optimal prediction problem is more precisely formulated in (2.2) below and the main result is given in Theorem 1.

It was recently observed in [11] that the optimal stopping time in (2.2) when $\mu = 0$ is equal to $\tau_* = \{ t \in [0,T] \mid |B_t| \geq z_* \sqrt{T-t} \}$ where $z_* > 0$ is a specified constant and $B = (B_t)_{0 \leq t \leq T}$ is a standard Brownian motion (with no drift). This stopping time is equal in law to the stopping time $\tilde{\tau}_* = \inf \{ t \in [0,T] \mid S_t - B_t \geq z_* \sqrt{T-t} \}$ which is known (cf. [4]) to minimise $\mathsf{E}(S_T - B_\tau)^2$ over all stopping times $\tau$ of $B$ taking values in $[0,T]$, where we set $S_t = \max_{0 \leq s \leq t} B_s$ for $t \in [0,T]$. Therefore when there is no drift, the space domain problem of stopping as close as possible to the ultimate maximum of Brownian motion, and the time domain problem of stopping as close as possible to the last zero of Brownian motion, are exactly equivalent. It is natural to consider whether this connection is preserved when the drift $\mu \neq 0$. In this case it is known (cf. [3]) that the solution to the space domain problem (of minimising $\mathsf{E}(S_T^\mu - B_\tau^\mu)^2$) is given by $\tilde{\tau}_*^\mu = \inf \{ t \in [0,T] \mid b_1(t) \leq S_t^\mu - B_t^\mu \leq b_2(t) \}$ where $b_1$ and $b_2$ are specified functions of time and $t_* \in [0,T)$. This stopping time has a distribution which is rather different from that of $\tilde{\tau}_*$ above (when $\mu$ is away from zero) raising the question of whether the solution to (2.2) below will bare any resemblance to it. In this paper we show that the optimal stopping time in (2.2) is in fact more similar to $\tau_*$ above. Predicting the last zero of Brownian motion with drift is therefore quite distinct from predicting its ultimate maximum.

Turning to a general stochastic process $X = (X_t)_{0 \leq t \leq T}$ one can ask similar questions about predicting the last zero of $X$ in the interval $[0,T]$. Certainly, if $X$ is discontinuous or multidimensional one would expect the problem to be much harder. Although these issues are not directly addressed in the paper, a simple approach based on the local time of $X$ is briefly presented in Remark 1 below. This may be used as a starting point from which to study similar problems for more general processes.

## 2. Formulation of the problem

We begin our exposition by formally introducing the setting and the problem to be studied. Let $B = (B_t)_{0 \leq t \leq T}$ be a standard Brownian motion defined on the probability space $(\Omega, \mathcal{F}, \mathsf{P})$ with $B_0 = 0$ under $\mathsf{P}$. For any $\mu \in I\!R$ let $B^\mu = (B_t^\mu)_{0 \leq t \leq T}$ be the Brownian motion with drift $\mu$ given by $B_t^\mu = \mu t + B_t$ for $t \in [0,1]$. Define the last zero of $B^\mu$ as

$$g = \sup \{ 0 \leq t \leq T \mid B_t^\mu = 0 \} \tag{2.1}$$

so that $g$ is well–defined under $\mathsf{P}$ since $B^\mu$ starts at zero under this measure. Later on we shall need to allow $B^\mu$ to start at arbitrary values in $I\!R$, raising the question of which value to assign to $g$ if $B^\mu$ does not hit zero at all on $[0,T]$. However, as it turns out (see Lemma 1 below), the value of $g$ in this case is irrelevant.

Consider the optimal prediction problem

$$V_* = \inf_{0 \leq \tau \leq T} \mathsf{E}|g - \tau| \tag{2.2}$$

where the infimum is taken over all stopping times $\tau$ of $B^\mu$. The gain process $(|g-t|)_{0 \leq t \leq T}$ is not adapted to the natural filtration $(\mathcal{F}_t^{B^\mu})_{0 \leq t \leq T}$ of $B^\mu$ as $g$ is only $\mathcal{F}_T^{B^\mu}$ measurable, meaning that the problem above falls outside the scope of standard optimal stopping theory. However, it is well known that there exist close links between $g$ and the maximum of $B^\mu$, and these we now exploit with the aim of reducing (2.2) to an equivalent optimisation problem to which the standard techniques of optimal stopping for Markov processes (see e.g. [7]) can be applied.

Let $S^\mu = (S_t^\mu)_{0 \le t \le T}$ be defined by $S_t^\mu = \max_{0 \le s \le t} B_s^\mu$ so that $S_t^\mu$ denotes the maximum of $B^\mu$ up to time $t \in [0,T]$. Then it follows (cf. [2] & [5]) that the distribution of $S_t^\mu$ is given by the formula

$$F^{(\mu)}(t,x) := \mathsf{P}(S_t^\mu \le x) = \Phi\Big(\frac{x-\mu t}{\sqrt{t}}\Big) - e^{2\mu x}\,\Phi\Big(\frac{-x-\mu t}{\sqrt{t}}\Big) \tag{2.3}$$

for all $(t,x) \in [0,\infty) \times I\!R_+$ where $\Phi(x) = \int_{-\infty}^x \varphi(z)\,dz$ denotes the distribution function of a standard normal random variable and $\varphi(x) = (1/\sqrt{2\pi})\,e^{-x^2/2}$ denotes its density for $x \in I\!R$.

**Lemma 2.1:** *The optimal prediction problem* (2.2) *above is equivalent to the standard optimal stopping problem*

$$V = \inf_{0 \le \tau \le T} \mathsf{E}\Big(\int_0^\tau H(t, B_t^\mu)\,dt\Big) \tag{2.4}$$

*where the function $H$ is given by*

$$H(t,x) = 2F^{(-\mu)}(T-t,x)\,I(x>0) + 2F^{(\mu)}(T-t,-x)\,I(x<0) - 1 \tag{2.5}$$

*for all $(t,x) \in [0,T] \times I\!R$, and where $V_*$ from* (2.2) *is given by $V_* = V + \mathsf{E}(g)$.*

**Proof:** We follow the same approach as in [11, Lemma 1] and [9]. Fix any stopping time $\tau$ of $B^\mu$ taking values in $[0,T]$. We then have

$$\begin{aligned} |g-\tau| &= (\tau-g)^+ + (\tau-g)^- = (\tau-g)^+ + g - \tau \wedge g \\ &= \int_0^\tau I(g \le t)\ dt + g - \int_0^\tau I(g > t)\ dt \\ &= g + \int_0^\tau \big[2\,I(g \le t) - 1\big]\,dt\,. \end{aligned} \tag{2.6}$$

From Fubini's theorem one finds that

$$\begin{aligned} \mathsf{E}\Big[\int_0^\tau I(g \le t)\ dt\Big] &= \int_0^\infty \mathsf{E}\,[\,I(\tau > t)\,I(g \le t)\,]\,dt \\ &= \mathsf{E}\Big[\int_0^\infty \mathsf{E}\,[\,I(\tau > t)\,I(g \le t)\,|\,\mathcal{F}_t^{B^\mu}]\,dt\Big] \\ &= \mathsf{E}\Big[\int_0^\tau \mathsf{P}(g \le t\,|\,\mathcal{F}_t^{B^\mu})\,dt\Big]\,. \end{aligned} \tag{2.7}$$

Using the continuity of the sample paths of $B^\mu$ to compute the conditional probability above we arrive at the following expression:

$$\mathsf{P}(g \le t\,|\,\mathcal{F}_t) = \begin{cases} \mathsf{P}\big(\min\limits_{t \le r \le T} B_r^\mu > 0\,\big|\,\mathcal{F}_t^{B^\mu}\big) & \text{if } \ B_t^\mu > 0 \\ \mathsf{P}\big(\max\limits_{t \le r \le T} B_r^\mu < 0\,\big|\,\mathcal{F}_t^{B^\mu}\big) & \text{if } \ B_t^\mu < 0 \end{cases} \tag{2.8}$$

for $t \in [0,T]$ (the case $B_t^\mu = 0$ can be ignored since the amount of time Brownian motion spends at zero has Lebesgue measure zero). If we consider the case when

$B_t^\mu > 0$ we find that

$$\begin{aligned} \mathsf{P}\big(\min_{t\le s\le T} B_s^\mu > 0 \,\big|\, \mathcal{F}_t^{B^\mu}\big) &= \mathsf{P}\big(\min_{t\le s\le T} (B_s^\mu - B_t^\mu) + B_t^\mu > 0 \,\big|\, \mathcal{F}_t^{B^\mu}\big) \\ &= \mathsf{P}\big(\max_{t\le s\le T} (-B_s^\mu + B_t^\mu) < b \,\big|\, \mathcal{F}_t^{B^\mu}\big)\Big|_{b=B_t^\mu} \\ &= \mathsf{P}\big(S_{T-t}^{-\mu} < b\big)\Big|_{b=B_t^\mu} \end{aligned} \tag{2.9}$$

for any $t \in [0,T]$ since $B^\mu$ has stationary independent increments. A similar argument when $B_t^\mu < 0$ shows that

$$\mathsf{P}\big(\max_{t\le s\le T} B_s^\mu < 0 \,\big|\, \mathcal{F}_t^{B^\mu}\big) = \mathsf{P}\big(S_{T-t}^{\mu} < -b\big)\Big|_{b=B_t^\mu} \tag{2.10}$$

for any $t \in [0,T]$. Taking expectations on both sides of (2.6) and using (2.7)–(2.10) we get

$$\begin{aligned} \mathsf{E}|g-\tau| = \mathsf{E}(g) + \mathsf{E}\bigg[\int_0^\tau \Big[&2F^{(-\mu)}(T\!-\!s, B_s^\mu)\, I(B_s^\mu > 0) \\ &+ 2F^{(\mu)}(T\!-\!s, -B_s^\mu)\, I(B_s^\mu < 0) - 1\Big]\, ds\bigg] \end{aligned} \tag{2.11}$$

for any stopping time $\tau$ of $B^\mu$ with values in $[0,T]$. Taking the infimum over all such stopping times we conclude the proof. □

In order to apply the standard techniques from the theory of optimal stopping for Markov processes (see e.g. [7]) it is necessary to extend the problem (2.4) by allowing $B^\mu$ to start at any time $t \in [0,T]$ at any point $x$ in the state space. Since $B^\mu$ is a time–homogeneous Markov process there exists a family of Markov measures $(\mathsf{P}_x)_{x\in I\!R}$ such that $\mathsf{P}_x(B_0^\mu = x) = 1$ for any $x \in I\!R$. Using these measures to change the starting point of the process results in the problem (2.4) becoming

$$V(t,x) = \inf_{0\le\tau\le T-t} \mathsf{E}_x\bigg[\int_0^\tau H(t\!+\!s, B_s^\mu)\, ds\bigg] \tag{2.12}$$

for any $(t,x) \in [0,T]\times I\!R$ so that the value $V_*$ of the optimal stopping problem (2.2) is given by $V_* = V(0,0) + \mathsf{E}(g)$. However, it is well known that the process $(x\!+\!B_t^\mu)_{t\ge 0}$ under $\mathsf{P}_0$ provides an explicit realisation of $B^\mu$ under $\mathsf{P}_x$ for any $x \in I\!R$. This together with the fact that the optimal stopping time in (2.4) is the first hitting time to a set (this follows from general theory and will also be proved below) leads to the further simplification

$$V(t,x) = \inf_{0\le\tau\le T-t} \mathsf{E}\bigg[\int_0^\tau H(t\!+\!s, x\!+\!B_s^\mu)\, ds\bigg] \tag{2.13}$$

where $\mathsf{P}$ is the original measure $\mathsf{P}_0$ under which $B^\mu$ starts at zero at time 0. This expression then is the form we shall use when solving the optimal prediction problem (2.2). Note that $V \le 0$ since we can always insert the stopping time identically equal to zero in the expectation above.

## 3. Result and proof

In order to state our main result we need the following definitions. Let the function $f$ denote the density of $B_t^\mu$ so that

$$f(t,x) = \frac{1}{\sqrt{2\pi t}}\, e^{-\frac{1}{2t}(x-\mu t)^2} \tag{3.1}$$

for all $(t,x) \in (0,T]\times I\!R$. Recalling $H$ in (2.5) above define the function $K$ as

$$\begin{aligned} K(t,x,s,z_-,z_+) &= \mathsf{E}\big[H(t+s,x+B_s^\mu)\, I(z_- < x+B_s^\mu < z_+)\big] \\ &= \int_{-x+z_-}^{-x+z_+} H(t+s,x+b)\, f(s,b)\, db \end{aligned} \tag{3.2}$$

for all $(t,x) \in [0,T]\times I\!R$, all $s \in [0,T-t]$ and $z_- < z_+$ in $I\!R$. Lastly the set $\{H<0\} := \{\,(t,x)\in[0,T]\times I\!R \mid H(t,x)<0\,\}$ will play a prominent role in our discussion. A direct exam- ination of $H$ reveals the existence of two continuous functions $h_-$ and $h_+$ mapping $[0,T]$ into $I\!R$ such that $\{H<0\} = \{\,(t,x)\in[0,T]\times I\!R \mid h_-(t)<x<h_+(t)\,\}$. Furthermore $h_-$ is increasing and $h_+$ is decreasing with $h_- \le 0 \le h_+$ and $h_-(T) = h_+(T) = 0$. We may now state our main result.

**Theorem 3.1:** *Consider the optimal stopping problem* (2.13). *Then there exist continuous functions $b_-$ and $b_+$ on $[0,T]$ such that the optimal stopping set is given by*

$$D = \{\,(t,x)\in[0,T]\times I\!R \mid x \le b_-(t) \;\; or \;\; x \ge b_+(t)\,\} \tag{3.3}$$

*where $t\mapsto b_-(t)$ is increasing, $t\mapsto b_+(t)$ is decreasing, and $b_-(T) = b_+(T) = 0$. This means that the stopping time*

$$\tau_D(t,x) = \inf\{\, s\in[0,T-t] \mid x+B_s^\mu \in D\,\} \tag{3.4}$$

*is optimal for all $(t,x)\in[0,T]\times I\!R$. The value function $V$ from* (2.13) *is given by*

$$V(t,x) = \int_0^{T-t} K\big(t,x,s,b_-(t+s),b_+(t+s)\big)\, ds \tag{3.5}$$

*for all $(t,x)\in[0,T]\times I\!R$, and the functions $b_-$ and $b_+$ themselves are uniquely characterised by the coupled system of non–linear Volterra integral equations*

$$\int_0^{T-t} K\big(t,b_-(t),s,b_-(t+s),b_+(t+s)\big)\, ds = 0 \tag{3.6}$$

$$\int_0^{T-t} K\big(t,b_+(t),s,b_-(t+s),b_+(t+s)\big)\, ds = 0 \tag{3.7}$$

*in the sense that they are the unique solution to* (3.6)+(3.7) *in the class of continuous functions $t\mapsto b_-(t)$ and $t\mapsto b_+(t)$ on $[0,T]$ satisfying $b_-(t)\le h_-(t)$ and $b_+(t)\ge h_+(t)$ for all $t\in[0,T]$. Finally, the value $V_*$ from* (2.2) *is given by $V_* = V(0,0)+\mathsf{E}(g)$ and the optimal stopping time for this problem is $\tau_D(0,0)$.*

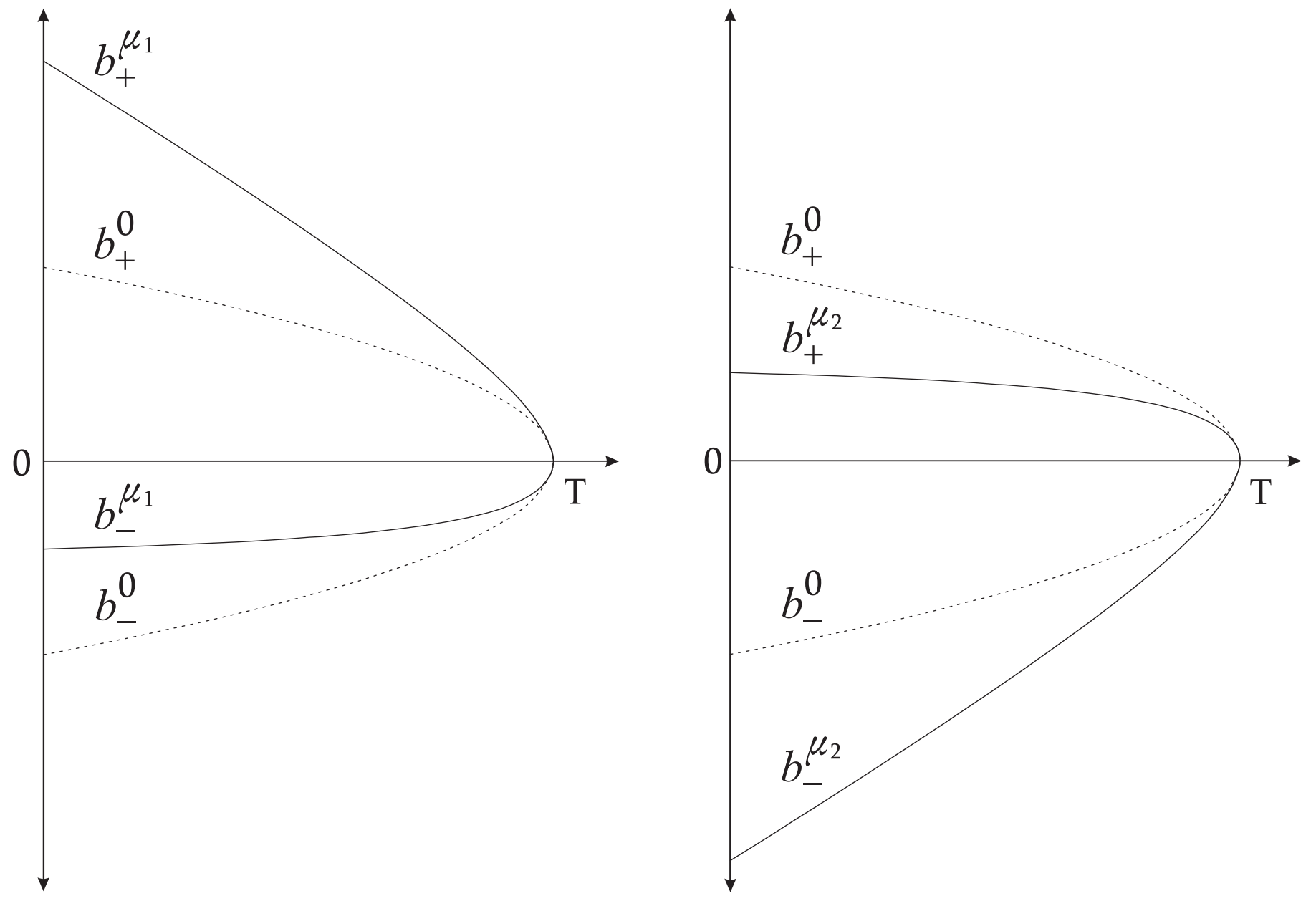

Figure 1. A computer drawing of the optimal stopping boundaries for Brownian motion with drift $\mu_1 < 0$ and $\mu_2 > 0$. The dotted line is the optimal stopping boundary for Brownian motion with zero drift.

**Proof:** 1. *Existence of optimal stopping time.* We begin by showing that an optimal stopping time for the problem (2.13) exists. Since $H$ is continuous and bounded, and the flow $x \mapsto x + B^\mu$ is continuous, it follows that the map $(t,x) \mapsto \mathsf{E}[\int_0^\tau H(t+r, x+B_r^\mu)\,dr]$ is continuous and thus upper semicontinuous (usc) for any stopping time $\tau$ with values in $[0,T]$. The function $V$ is therefore usc as well (being the infimum of usc functions) and so by general results of optimal stopping (see [7, Corollary 2.9]) it follows that an optimal stopping time for the problem (2.13) exists. Moreover, this stopping time is given by (3.4) above where the stopping set equals $D = \{\,(t,x) \in [0,T]\times I\!R \mid V(t,x) = 0\,\}$ and the continuation set is given by $C = \{\,(t,x) \in [0,T]\times I\!R \mid V(t,x) < 0\,\} = D^c$. The fact that $D$ is closed (and $C$ is open) follows from the fact that $V$ is usc.

2. *Shape of $D$.* From the functional form of $H$ given at (2.5) (see also (2.3) above), one sees that the map $t \mapsto H(t,x)$ is increasing on $[0,T]$ for every $x \in I\!R$. Fix any $x \in I\!R$ and $s < t$ in $[0,T]$ and set $\tau_s = \tau_D(s,x)$ and $\tau_t = \tau_D(t,x)$. Since $0 \le \tau_t \le T-t < T-s$ one finds that

$$\begin{aligned} &V(t,x) - V(s,x) \qquad (3.8)\\ &= \mathsf{E}\left[\int_0^{\tau_t} H(t+r, x+B_r^\mu)\,dr\right] - \mathsf{E}\left[\int_0^{\tau_s} H(s+r, x+B_r^\mu)\,dr\right]\\ &\ge \mathsf{E}\left[\int_0^{\tau_t} H(t+r, x+B_r^\mu) - H(s+r, x+B_r^\mu)\,dr\right] \ge 0 \end{aligned}$$

from where we derive the important fact that

$$t \mapsto V(t,x) \text{ is increasing on } [0,T] \tag{3.9}$$

for each $x \in I\!R$ given and fixed. A direct consequence of this is that if any point $(t,x) \in D$, then all the points $(t+s,x) \in D$ for $s \in [0,T-t]$ since $0 \ge V(t+s,x) \ge V(t,x) = 0$.

Turning to the continuation set $C$, take any point $(t,x)$ in the open set $\{H<0\}$ and let $U \subset \{H<0\}$ be any open neighbourhood of $(t,x)$. Defining the stopping time $\sigma_U$ as the first exit time from $U$, we see that

$$V(t,x) \le \mathsf{E}\left[\int_0^{\sigma_U} H(t+s, x+B_s^\mu)\, ds\right] < 0 \tag{3.10}$$

showing that $\{H<0\} \subseteq C$. We now show that $D$ has the form given in (3.3). Fix any $t$ in $[0,T)$ and $x \ge 0$ and suppose that $(t,x) \in D$. For any $y > x$ consider the process $B^\mu$ started at $(t,y)$. Since $(s,x) \in D$ for all $t \le s \le T$ it follows that $B^\mu$ started at $(t,y)$ must enter $D$ before crossing the level $x$, in other words $y+B_s^\mu \ge x$ for all $0 \le s \le \tau_D(t,y)$. From equation (2.5) we see that the map $x \mapsto H(t,x)$ is increasing on $I\!R_+$ for each $t \in [0,T]$ given and fixed. Hence for $y \ge x \ge 0$ we have

$$\begin{aligned} V(t,y) &= \mathsf{E}\left[\int_0^{\tau_D(t,y)} H(t+s, y+B_s^\mu)\, ds\right] \\ &\ge \mathsf{E}\left[\int_0^{\tau_D(t,y)} H(t+s, x+B_s^\mu)\, ds\right] \ge \mathsf{E}\left[\int_0^{\tau_D(t,x)} H(t+s, x+B_s^\mu)\, ds\right] = 0 \end{aligned} \tag{3.11}$$

since $\tau_D(t,x) \equiv 0$. Therefore if $(t,x) \in D$ with $x \ge 0$, then $(s,y) \in D$ for all $y \ge x$ and $s \in [t,T]$. For $x \le 0$ a similar argument can be made showing that if $(t,x) \in D$ then $(s,y) \in D$ for all $y \le x$ and $s \in [t,T]$. It follows then that $D$ has the form given in (3.3) where the functions $b_-$ and $b_+$ are defined as

$$\begin{aligned} &b_+(t) = \inf\{\, x \ge 0 \mid (t,x) \in D \,\} \quad \text{and} \\ &b_-(t) = \sup\{\, x \le 0 \mid (t,x) \in D \,\} \end{aligned} \tag{3.12}$$

for all $t \in [0,T]$. We also see that $b_+$ is decreasing and $b_-$ is increasing.

3. *Functions $b_+$ and $b_-$ are finite valued.* We show that the optimal stopping boundaries $b_+$ and $b_-$ are finite valued, which also shows that the stopping region $D$ is strictly greater than the set $\{\,(T,x) \mid x \in I\!R\,\}$. We shall only present the proof for $b_+$ since the argument for $b_-$ is analogous.

Suppose that the function $b_+(t)$ is not finite valued for all $t \in [0,T]$ and define the time $t_* \in [0,T]$ as $t_* = \sup\{\, t \in [0,T] \mid b_+(t) = \infty \,\}$. There are two possibilities: either $b_+$ has a jump discontinuity at $t_*$ jumping down from infinity to a finite value (including the case when $t_* = T$), or $b_+$ has an asymptote at $t_*$. Suppose first that $b_+$ jumps from infinity at $t_*$ and that $t_* \in (0,T)$. If $\tau_x = \tau_D(0,x)$ and $\sigma_x = \inf\{\, t \in [0,T] \mid x+B_t^\mu < h_+(0)+1 \,\}$, then it is clear that $\tau_x \to t_*$ and $\sigma_x \to T$ as $x \to \infty$. Therefore we see that

$$\begin{aligned} 0 \ge V(0,x) &= \mathsf{E}\left[\int_0^{\tau_x} H(s, x+B_s^\mu)\, ds \; I(\tau_x \le \sigma_x)\right] \\ &\quad + \mathsf{E}\left[\int_0^{\tau_x} H(s, x+B_s^\mu)\, ds \; I(\tau_x > \sigma_x)\right] \\ &\ge \mathsf{E}\big[\varepsilon\, \tau_x\, I(\tau_x \le \sigma_x)\big] - T\,\mathsf{E}\big[\, I(\tau_x > \sigma_x)\big] \end{aligned} \tag{3.13}$$

where $\varepsilon := \inf\{\, H(t,x) \mid (t,x) \in [0,T] \times [h_+(0)+1, \infty) \,\} > 0$ and where we have used the fact that $H \ge -1$. Letting $x \to \infty$ and using the dominated convergence

theorem we obtain

$$0 \geq \lim_{x \to \infty} V(0,x) \geq \varepsilon\, t_* > 0 \tag{3.14}$$

which is a contradiction. This shows that we cannot have $t_* \in (0,T)$. However, since $b_+$ is right-continuous (this will be shown below) we cannot have $t_* = 0$ either, and so $b_+$ cannot have a jump discontinuity at $t_* \in [0,T)$. Finally, if $t_* = T$ then letting $x \to \infty$ and using the dominated convergence theorem we find that

$$\begin{aligned} 0 \geq V(0,x) &= \lim_{x \to \infty} \mathsf{E}\Big[ \int_0^{\tau_x} H(s, x{+}B_s^\mu)\, ds \Big] \\ &= \int_0^T H(s,\infty)\, ds = T > 0 \end{aligned} \tag{3.15}$$

which is a contradiction. This excludes the case $t_* = T$ and establishes the claim.

Suppose now that $b_+$ asymptotes to infinity at $t_*$ and that $t_* \in (0,T)$. It then follows that $t_* \leq \liminf_{x \to \infty} \tau_x \leq \limsup_{x \to \infty} \tau_x \leq t_* + \delta < T$ for some $\delta > 0$, and the same argument as in (3.13) yields $0 \geq \limsup_{x \to \infty} V(0,x) \geq \varepsilon\, t_* > 0$. Lastly, suppose that $b_+$ asymptotes to infinity at 0. By extending the terminal time to $T' > T$ and considering our problem on the interval $[0,T']$ instead of $[0,T]$ we will (since $b_+$ is decreasing) shift $t_*$ to be strictly positive, reducing it to the case already considered. Therefore $b_+(t)$ must be finite for all $t \in [0,T]$.

4. *Continuity of V*. We show that $(t,x) \mapsto V(t,x)$ is continuous on $[0,T] \times I\!R$. For this, take any $t \in [0,T]$ and $x,y \in I\!R$, set $\tau_x = \tau_D(t,x)$ and $\tau_y = \tau_D(t,y)$ and suppose without loss of generality that $x \leq y$. It follows then that

$$\begin{aligned} &\mathsf{E}\Big[ \int_0^{\tau_y} H(t{+}s, y{+}B_s^\mu) - H(t{+}s, x{+}B_s^\mu)\, ds \Big] \\ &\leq V(t,y) - V(t,x) \leq \mathsf{E}\Big[ \int_0^{\tau_x} H(t{+}s, y{+}B_s^\mu) - H(t{+}s, x{+}B_s^\mu)\, ds \Big]. \end{aligned} \tag{3.16}$$

Turning to equation (2.5) and manipulating the quantities we find that

$$\begin{aligned} &H(t,y) - H(t,x) \\ &= 2\Big( F^{(-\mu)}(T{-}t, y) - F^{(-\mu)}(T{-}t, x) \Big)\, I(x > 0) \\ &- 2\Big( F^{(\mu)}(T{-}t, -x) - F^{(\mu)}(T{-}t, -y) \Big)\, I(y < 0) \\ &+ 2\Big( F^{(-\mu)}(T{-}t, y) - F^{(\mu)}(T{-}t, -x) \Big)\, I(x \leq 0 \leq y)\, . \end{aligned} \tag{3.17}$$

Applying the mean value theorem in the first and second line upon noting that

$$F_x^{(\mu)}(t,x) = \frac{2}{\sqrt{t}}\, \varphi\Big( \frac{x - \mu t}{\sqrt{t}} \Big) - 2\mu e^{2\mu x}\, \Phi\Big( \frac{-x - \mu t}{\sqrt{t}} \Big) \leq \frac{2}{\sqrt{t}} + 2\,|\mu| \tag{3.18}$$

for all $(t,x) \in (0,T] \times I\!R_+$ we obtain the inequality

$$
\begin{aligned}
&-4(y-x)\Big(\frac{1}{\sqrt{T-t}} + |\mu|\Big)\, I(y<0) \;-\; 2\, I(x \le 0 \le y) && (3.19)\\
&\le H(t,y) - H(t,x)\\
&\le 4(y-x)\Big(\frac{1}{\sqrt{T-t}} + |\mu|\Big)\, I(x>0) \;+\; 2\, I(x \le 0 \le y)
\end{aligned}
$$

which we can insert into equation (3.16) to get

$$
\begin{aligned}
&|V(t,y) - V(t,x)| && (3.20)\\
&\le \mathsf{E}\Big[\int_0^{T-t} \Big(4(y-x)\big((T-t-s)^{-1/2} + |\mu|\big)\\
&\qquad\qquad + 2\, I(x+B_s^\mu \le 0 \le y+B_s^\mu)\Big)\, ds\Big]\\
&= 4(y-x)\Big(2\sqrt{T-t} + |\mu|(T-t)\Big)\\
&\qquad\qquad + 2\int_0^{T-t} \Big(\Phi\Big(\frac{y+\mu s}{\sqrt{s}}\Big) - \Phi\Big(\frac{x+\mu s}{\sqrt{s}}\Big)\Big)\, ds\\
&\le 4(y-x)\Big(2\sqrt{T-t} + |\mu|(T-t)\Big) + 4\sqrt{T-t}\,(y-x)\\
&\le \Big(12\sqrt{T} + 4T|\mu|\Big)(y-x)
\end{aligned}
$$

upon using that $\Phi(z) - \Phi(w) \le z - w$ for $w \le z$ in $I\!R$. Hence we conclude that $x \mapsto V(t,x)$ is continuous on $I\!R$ uniformly over all $t \in [0,T]$.

It remains to show that $t \mapsto V(t,x)$ is continuous on $[0,T]$ for every $x \in I\!R$. To see this fix $x \in I\!R$, take any $s \le t$ in $[0,T]$ and set $\tau_s = \tau_D(s,x)$ and $\tau_t = \tau_D(t,x)$. Using the stopping time $\sigma = \tau_s \wedge (T-t)$ instead of $\tau_D(t,x)$ one finds that

$$
\begin{aligned}
&0 \le V(t,x) - V(s,x) && (3.21)\\
&\quad\le \mathsf{E}\Big[\int_0^\sigma H(t+r, x+B_r^\mu) - H(s+r, x+B_r^\mu)\, dr\Big]\\
&\qquad - \mathsf{E}\Big[\int_\sigma^{\tau_s} H(s+r, x+B^\mu)\, dr\Big].
\end{aligned}
$$

Taking the limit as $t - s \to 0$ and using the dominated convergence theorem, the continuity of $t \mapsto H(t,x)$ and the fact that $0 \le \tau_s - \sigma \le t - s$, we see that both expectations on the right hand side vanish and we conclude that $(t,x) \mapsto V(t,x)$ is continuous on $[0,T] \times I\!R$.

5. *Free-boundary problem.* We now formulate a free-boundary problem that the value function $V$ solves. This differential equation will be useful to us later on, but is also interesting in its own right and can be used as the departure point in computing numerical values for the boundaries $b_+$ and $b_-$ and for the value function $V$. It is well known from the theory of Markov processes (see e.g. [7, Chapter III, Section 7]) that $V$ is $C^{1,2}$ in the continuation region $C$ and satisfies the following version of the *Kolmogorov backward equation*:

$$V_t(t,x) + \mu V_x(t,x) + \tfrac{1}{2}V_{xx}(t,x) = -H(t,x) \quad \text{for all } (t,x) \in C \tag{3.22}$$

$$V(t,x) = 0 \quad \text{for all } (t,x) \in D. \tag{3.23}$$

The last condition which forms part of the free-boundary problem is the smooth fit property which reads

$$x \mapsto V_x(t,x) \text{ is continuous over } b_-(t) \text{ and } b_+(t) \text{ for all } t \in [0,T). \tag{3.24}$$

This result, which we establish below, will play an important role in the derivation of the integral equations (3.5) and (3.6)+(3.7).

6. *Smooth fit.* We show that $x \mapsto V_x(t,x)$ exists and is continuous over the optimal stopping boundaries $b_-$ and $b_+$. Since the arguments for the two curves are very similar we will only present the proof for $b_+$.

Fix any time $t \in [0,T)$ and $\varepsilon > 0$, and set $x = b_+(t)$ and $\tau_\varepsilon = \tau_D(t, x-\varepsilon)$. From equations (3.16) and (3.19) we find that

$$\begin{aligned} 0 &\le V(t,x) - V(t,x-\varepsilon) \\ &\le \mathsf{E}\Big[\int_0^{\tau_\varepsilon} H(t,x+B_s^\mu) - H(t,x-\varepsilon+B_s^\mu)\,ds\Big] \\ &\le \mathsf{E}\Big[\int_0^{\tau_\varepsilon} \Big(4\varepsilon\Big((T-t-s)^{-1/2} + |\mu|\Big) + 2\,I(0 \le x+B_s^\mu \le \varepsilon)\Big)\,ds\Big] \\ &\le 4\varepsilon\,\mathsf{E}\Big[2\big(\sqrt{T-t} - \sqrt{T-t-\tau_\varepsilon}\,\big) + |\mu|\,\tau_\varepsilon\Big] \\ &\quad + 2\,\mathsf{E}\Big[\int_0^{\tau_\varepsilon} I(0 \le x+B_s^\mu \le \varepsilon)\,ds\Big]. \end{aligned} \tag{3.25}$$

Consider the final expectation above and denote by $\ell^x(B^\mu) = \big(\ell_t^x(B^\mu)\big)_{t\ge0}$ the local time of $B^\mu$ at the level $x$. Then it follows from the occupation times formula (see e.g. [8, p. 224]) and the integral mean value theorem that

$$\begin{aligned} \frac{1}{\varepsilon}\int_0^{\tau_\varepsilon} I(0 \le x+B_s^\mu \le \varepsilon)\,ds &= \frac{1}{\varepsilon}\int_{-\infty}^{\infty} I(-x \le a \le -x+\varepsilon)\,\ell_{\tau_\varepsilon}^a(B^\mu)\,da \\ &= \frac{1}{\varepsilon}\int_0^{\varepsilon} \ell_{\tau_\varepsilon}^{-x+a}(B^\mu)\,da = \ell_{\tau_\varepsilon}^{-x+\eta_\varepsilon}(B^\mu) \end{aligned} \tag{3.26}$$

where $\eta_\varepsilon$ is a random variable taking values in $[0,\varepsilon]$. Define the stopping time $\sigma_\varepsilon$ to be the first time $B^\mu$ rises above the level $x$ when $B^\mu$ starts at time $t$ at the point $x-\varepsilon$, in other words $\sigma_\varepsilon = \inf\{\, s \in [0,T-t] \mid x-\varepsilon+B_s^\mu = x \,\}$. Since the map $t \mapsto b_+(t)$ is decreasing it follows that $0 \le \tau_\varepsilon \le \sigma_\varepsilon$. Recalling that $t \mapsto -\mu t$ is a lower function of Brownian motion at $0+$ we see that $\sigma_\varepsilon$, and consequently $\tau_\varepsilon$, tend to zero as $\varepsilon \downarrow 0$.

Dividing equation (3.25) through by $\varepsilon$ and passing to the limit as $\varepsilon \downarrow 0$, we see that the first expectation disappears (upon using the dominated convergence theorem) so that

$$0 \le \lim_{\varepsilon\to0} \frac{V(t,x) - V(t,x-\varepsilon)}{\varepsilon} \le 2\lim_{\varepsilon\to0} \mathsf{E}\Big[\ell_{\tau_\varepsilon}^{-x+\eta_\varepsilon}(B^\mu)\Big] = 0 \tag{3.27}$$

where for the final equality we use the fact that $\ell^{-x+\eta_\varepsilon}_{\tau_\varepsilon}(B^\mu)$ is bounded by an integrable random variable when $\varepsilon > 0$ and the well-known result that $(t,a) \mapsto \ell^a_t(B^\mu)$ is continuous on $I\!R_+ \times I\!R$. The former can be seen by recalling the Tanaka formula

$$\begin{aligned} |B^\mu_t - a| &= |a| + \int_0^t \mathrm{sign}(B^\mu_s - a)\, dB^\mu_s + \ell^a_t(B^\mu) \\ &= |a| + \beta_t + \mu \int_0^t \mathrm{sign}(B^\mu_s - a)\, ds + \ell^a_t(B^\mu) \end{aligned} \tag{3.28}$$

where $a \in I\!R$ is any number and $\beta_t = \int_0^t \mathrm{sign}(B^\mu_s - a)\, dB_s$ for $t \geq 0$ is a standard Brownian motion by Lévy's characterisation theorem. It follows from (3.28) that

$$\sup_{0 \leq t \leq T} \sup_{a \in I\!R} \ell^a_t(B^\mu) \leq \max_{0 \leq t \leq T} |B_t| + \max_{0 \leq t \leq T} |\beta_t| + 2|\mu| T \tag{3.29}$$

and the claim above follows by noting that the right-hand side in (3.29) defines an integrable random variable. We therefore conclude that $x \mapsto V(t,x)$ is differentiable at $b_+(t)$ and that $V_x\big(t, b_+(t)\big) = 0$.

A small modification of the argument above shows that $x \mapsto V_x(t,x)$ is continuously differentiable at $b_+(t)$. Indeed, taking $\delta > 0$ and setting $\tau_\delta = \tau_D(t, x-\delta)$ where $x = b_+(t)$, we find using equations (3.16), (3.19) and (3.26) that

$$\begin{aligned} \frac{V(t, x-\delta+\varepsilon) - V(t, x-\delta)}{\varepsilon} &\leq 4\, \mathsf{E}\Big[2\big(\sqrt{T-t} - \sqrt{T-t-\tau_\delta}\,\big) + |\mu|\, \tau_\delta\Big] \\ &\quad + 2\, \mathsf{E}\Big[\ell^{-x+\delta-\eta_\varepsilon}_{\tau_\delta}(B^\mu)\Big] \end{aligned} \tag{3.30}$$

where $\varepsilon \in (0, \delta)$ is a given number and $\eta_\varepsilon$ is a random variable taking values in $[0, \varepsilon]$. In exactly the same way we find that

$$\begin{aligned} &- 4\, \mathsf{E}\Big[2\big(\sqrt{T-t} - \sqrt{T-t-\tau_\delta}\,\big) + |\mu|\, \tau_\delta\Big] - 2\, \mathsf{E}\Big[\ell^{-x+\delta+\rho_\varepsilon}_{\tau_\delta}(B^\mu)\Big] \\ &\leq \frac{V(t, x-\delta) - V(t, x-\delta-\varepsilon)}{\varepsilon} \end{aligned} \tag{3.31}$$

where $\varepsilon > 0$ is a given number and $\rho_\varepsilon$ is a random variable taking values in $[0, \varepsilon]$. Letting first $\varepsilon \downarrow 0$ in (3.30) and (3.31) (upon recalling that $V$ is $C^{1,2}$ in $C$ so that $V_x(t, x-\delta)$ exists) and then $\delta \downarrow 0$ in the resulting inequalities (upon recalling that $\tau_\delta \to 0$) we see that the map $x \mapsto V_x(t,x)$ is continuous at $b_+(t)$.

7. *Continuity of boundaries.* We show that the functions $t \mapsto b_-(t)$ and $t \mapsto b_+(t)$ are continuous on $[0,T]$. As before we will restrict ourselves to $b_+$ since the argument for $b_-$ is similar.

We first show that $b_+$ is right-continuous. Fix any $t \in [0,T)$, let $t_n \downarrow t$ and consider the limit $b_+(t+) := \lim_{n\to\infty} b_+(t_n)$ which exists as $b_+$ is decreasing. Since $\big(t_n, b_+(t_n)\big) \in D$ for all $n \geq 1$ and $D$ is closed, it follows that $\big(t, b_+(t+)\big) \in D$ and so from (3.12) we see that $b_+(t+) \geq b_+(t)$. On the other hand, the fact that $b_+$ is decreasing implies that $b_+(t) \geq b_+(t_n)$ for all $n \geq 1$, and passing to the limit as $n \to \infty$ we obtain the reverse inequality.

We now show that $b_+$ is left continuous. For this, suppose there exists some $t \in (0,T]$ at which $b_+(t-) > b_+(t)$ and choose any $x \in \big(b_+(t), b_+(t-)\big)$. Since $b_+ \geq h_+$ and $h_+$ is continuous, it follows that $x \geq h_+(s)$ for all $s \in [s_1, t]$ for some $s_1$ sufficiently close to $t$. Hence $m := \inf \big\{ H(s,y) \;\big|\; s \in [s_1,t),\ y \in [x, b_+(s)] \big\} > 0$ by the continuity of $H$. Moreover, since $V$ is continuous and $V(t,y) = 0$ for all $y \in [x, b_+(t-)]$, it follows that

$$|\mu V(s,y)| \le \frac{m}{4}\big[b_+(t-) - x\big] \tag{3.32}$$

for all $s \in [s_2, t]$ and $y \in [x, b_+(s)]$ where $s_2 \in [s_1, t)$ is some value sufficiently close to $t$. Therefore for any $s \in [s_2, t)$ we find using (3.22) and (3.24) that

$$\begin{aligned} V(s,x) &= \int_x^{b_+(s)} \int_y^{b_+(s)} V_{xx}(s,z)\, dz\, dy \\ &= -2 \int_x^{b_+(s)} \int_y^{b_+(s)} \big(V_t + \mu V_x + H\big)(s,z)\, dz\, dy \\ &\le 2 \int_x^{b_+(s)} \Big(\mu V(s,y) - m\big(b_+(s) - y\big)\Big)\, dy \\ &\le \frac{m}{2}\big(b_+(t-) - x\big)\big(b_+(s) - x\big) - m\big(b_+(s) - x\big)^2 \end{aligned} \tag{3.33}$$

where in the second last inequality we used $V_t \ge 0$ and in the last we used (3.32). Passing to the limit as $s \uparrow t$ we see that $V(t,x) \le -\frac{m}{2}\big(b_+(t-) - x\big)^2 < 0$ which contradicts the fact that $(t,x) \in D$. We conclude that $t \mapsto b_+(t)$ is continuous on $[0,T]$. Note that this proof also shows that $b_+(T) = 0 = b_-(T)$ since $h_+(T) = 0 = h_-(T)$ and $V(x,T) = 0$ for all $x \in I\!R$.

8. *Integral equations.* We may now derive the integral equations (3.5) and (3.6)+(3.7). From equation (3.30) we see that

$$|V_x(t,x)| \le 4(2\sqrt{T} + |\mu| T) + 2\,\mathsf{E}\big[Z\big] =: K < \infty \tag{3.34}$$

for all $(t,x) \in [0,T] \times I\!R$ where $Z := \sup_{a \in I\!R} \ell_T^a(B^\mu)$. Inserting this into (3.22) and using (3.9) we obtain $V_{xx} = -2(V_t + \mu V_x + H) \le 2(|\mu|\, K - H)$ in $C$. If we let

$$f(t,x) = 2 \int_0^x \int_0^y \big(1 + |\mu|\, K - H(t,z)\big)\, dz\, dy \tag{3.35}$$

for all $(t,x) \in [0,T] \times I\!R$ then we have $V_{xx} \le f_{xx}$ on $C \cup D^o$ since $|H| \le 1$. Defining therefore the function $F : [0,T] \times I\!R \to I\!R$ by $F(t,x) = V(t,x) - f(t,x)$, we see that: (i) the map $x \mapsto F(t,x)$ is concave on each of the intervals $(-\infty, b_-(t))$, $(b_-(t), b_+(t))$ and $(b_+(t), \infty)$ for every $t \in [0,T]$; (ii) the function $F$ is $C^{1,2}$ on $C \cup D^o$; (iii) the function $F_t + \mu F_x + \frac{1}{2} F_{xx}$ is locally bounded on $C \cup D^o$; and (iv) the maps $t \mapsto F_x\big(t, b_-(t) \pm\big) = -f_x\big(t, b_-(t)\big)$ and $t \mapsto F_x\big(t, b_+(t) \pm\big) = -f_x\big(t, b_+(t)\big)$ are continuous on $[0,T]$. Since the boundaries $b_-$ and $b_+$ are monotone and consequently of bounded variation, we may apply the local time–space formula [6] to $F(t+s, B^\mu_{t+s})$ and Itô's formula to $f(t+s, B^\mu_{t+s})$ since $f$ is $C^{1,2}$. Adding these two expressions and using (3.22) and (3.24) we obtain

$$
\begin{aligned}
V(t+s, B^{\mu}_{t+s}) &= V(t,x) \\
&+ \int_0^s \big(V_t + \mu V_x + \tfrac{1}{2} V_{xx}\big)(t+r, B^{\mu}_{t+r}) \times \\
&\qquad\qquad I\big(B^{\mu}_{t+r} \notin \{b_-(t+r), b_+(t+r)\}\big)\; dr \\
&+ \int_0^s V_x(t+r, B^{\mu}_{t+r})\, I\big(B^{\mu}_{t+r} \notin \{b_-(t+r), b_+(t+r)\}\big)\; dB_{t+r} \\
&+ \frac{1}{2} \int_0^s \Big(V_x(t+r, B^{\mu}_{t+r}+) - V_x(t+r, B^{\mu}_{t+r}-)\Big) \times \\
&\qquad\qquad I\big(B^{\mu}_{t+r} = b_-(t+r)\big)\; d\ell^{b_-}_{t+r}(B^{\mu}) \\
&+ \frac{1}{2} \int_0^s \Big(V_x(t+r, B^{\mu}_{t+r}+) - V_x(t+r, B^{\mu}_{t+r}-)\Big) \times \\
&\qquad\qquad I\big(B^{\mu}_{t+r} = b_+(t+r)\big)\; d\ell^{b_+}_{t+r}(B^{\mu}) \\
&= V(t,x) - \int_0^s H(t+r, B^{\mu}_{t+r})\, I\big(b_-(t+r) < B^{\mu}_{t+r} < b_+(t+r)\big)\; dr + M_s
\end{aligned}
\tag{3.36}
$$

under $\mathsf{P}_{t,x}$ for any $(t,x) \in [0,T] \times I\!R$ and $s \in [0, T-t]$, where $\ell^{b_-}(B^{\mu})$ and $\ell^{b_+}(B^{\mu})$ are the local times of $B^{\mu}$ on the curves $b_-$ and $b_+$ respectively, and $M_s = \int_0^s V_x(t+r, B^{\mu}_{t+r})\, dB_{t+r}$ is a continuous martingale under $\mathsf{P}_{t,x}$ for $s \in [0, T-t]$. Setting $s = T-t$, taking $\mathsf{E}_{t,x}$ on both sides and using the optional sampling theorem we obtain

$$
V(t,x) = \mathsf{E}_{t,x}\bigg[\int_0^{T-t} H(t+r, B^{\mu}_{t+r})\, I\big(b_-(t+r) < B^{\mu}_{t+r} < b_+(t+r)\big)\; dr\bigg] \tag{3.37}
$$

which is exactly (3.5) after interchanging the order of integration. Setting $x$ equal to $b_-(t)$ and $b_+(t)$ in (3.37) we get

$$
\int_0^{T-t} \mathsf{E}_{t,b_-(t)}\Big[H(t+r, B^{\mu}_{t+r})\, I\big(b_-(t+r) < B^{\mu}_{t+r} < b_+(t+r)\big)\Big]\; dr = 0 \tag{3.38}
$$

$$
\int_0^{T-t} \mathsf{E}_{t,b_+(t)}\Big[H(t+r, B^{\mu}_{t+r})\, I\big(b_-(t+r) < B^{\mu}_{t+r} < b_+(t+r)\big)\Big]\; dr = 0 \tag{3.39}
$$

which is exactly (3.6)+(3.7) as claimed.

9. *Uniqueness*. Lastly we show that $b_-$ and $b_+$ are the unique solutions to (3.6)+(3.7) in the class of continuous functions $t \mapsto b_-(t)$ and $t \mapsto b_+(t)$ on $[0,T]$ satisfying $b_-(t) \le h_-(t)$ and $b_+(t) \ge h_+(t)$ for all $t$ in $[0,T]$.

Take any two continuous functions $c_+$ and $c_-$ on $[0,T]$ which solve (3.6)+(3.7) and satisfy $c_-(t) \le h_-(t)$ and $h_+(t) \le c_+(t)$ for all $t \in [0,T]$. Motivated by equation (3.37) above, define the continuous function $U^c : [0,T] \times I\!R \to I\!R$ as

$$
U^c(t,x) = \mathsf{E}_{t,x}\bigg[\int_0^{T-t} H(t+r, B^{\mu}_r)\, I\big(c_-(t+r) < B^{\mu}_{t+r} < c_+(t+r)\big)\; dr\bigg] \tag{3.40}
$$

and observe that $c_-$ and $c_+$ solving (3.6)+(3.7) means exactly that $U^c\big(t, c_-(t)\big) = U^c\big(t, c_+(t)\big) = 0$ for all $t \in [0,T]$. Let $D_c := \{\, (t,x) \in [0,T] \times I\!R \mid x \le c_-(t) \text{ or } x \ge c_+(t)\}$ so that $D_c$ is closed and plays the role of a 'stopping region' for $c_-$ and $c_+$. To avoid confusion we will denote by $D_b$ the original stopping set from (3.3) defined by the functions $b_-$ and $b_+$.

(i) We show that $U^c = 0$ on $D_c$. Since $B^\mu$ is Markov, the process

$$N_s := U^c(t{+}s, B^\mu_{t+s}) + \int_0^s H(t{+}r, B^\mu_{t+r})\, I\big(c_-(t{+}r) < B^\mu_{t+r} < c_+(t{+}r)\big)\, dr \tag{3.41}$$

is a martingale under $\mathsf{P}_{t,x}$ for all $s \in [0, T{-}t]$ and $x \in I\!R$. Take any point $(t,x) \in D_c$ and consider the stopping time

$$\sigma_c = \inf\big\{\, s \in [0, T{-}t] \;\big|\; B^\mu_{t+s} \notin D_c \,\big\} \tag{3.42}$$

under the measure $\mathsf{P}_{t,x}$. Since $U^c$ is zero on the curves $c_+$ and $c_-$ and $U^c(T,x) = 0$ for all $x \in I\!R$, we must have $U^c(t{+}\sigma_c, B^\mu_{t+\sigma_c}) = 0$. Inserting $\sigma_c$ in (3.41), taking $\mathsf{P}_{t,x}$ expectations and using the optional sampling theorem (since $H$ is bounded) we find that

$$U^c(t,x) = \mathsf{E}_{t,x}\Big[U^c(t{+}\sigma_x, B^\mu_{t+\sigma_c})\Big] = 0 \tag{3.43}$$

showing that $U^c = 0$ on $D_c$ as claimed.

(ii) We show that $U^c(t,x) \ge V(t,x)$ for all $(t,x) \in [0,T] \times I\!R$. To see this take any $(t,x) \in [0,T] \times I\!R$ and consider the stopping time

$$\tau_c = \inf\big\{\, s \in [0, T{-}t] \;\big|\; B^\mu_{t+s} \in D_c \,\big\} \tag{3.44}$$

under $\mathsf{P}_{t,x}$. We claim that $U^c(t{+}\tau_c, B^\mu_{t+\tau_c}) = 0$. Indeed, if $(t,x) \in D_c$ then $\tau_c = 0$ so that $U^c(t,x) = 0$ by the argument above. Conversely if $(t,x) \notin D_c$, then the result follows since $U^c$ is zero on $c_-$ and $c_+$ and $U^c(T,x) = 0$ for all $x \in I\!R$. Inserting $\tau_c$ in (3.41) and using the optional sampling theorem, we see that

$$\begin{aligned} U^c(t,x) &= \mathsf{E}_{t,x}\Big[\int_0^{\tau_c} H(t{+}s, B^\mu_{t+s})\, I\big(B^\mu_{t+s} \notin D_c\big)\; ds\Big] \\ &= \mathsf{E}_{t,x}\Big[\int_0^{\tau_c} H(t{+}s, B^\mu_{t+s})\, ds\Big] \ge V(t,x) \end{aligned} \tag{3.45}$$

where the second identity follows from the definition of $\tau_c$. We conclude that $U^c \ge V$ on $[0,T] \times I\!R$ as claimed.

(iii) We show that $D_b \subseteq D_c$. Suppose this is not the case so that there exists some time $t \in [0,T)$ for which either $b_+(t) < c_+(t)$ or $b_-(t) > c_-(t)$. Choose any $x > c_+(t)$ or $x < c_-(t)$ and consider the stopping time

$$\sigma_b = \inf\big\{\, s \in [0, T{-}t] \;\big|\; B^\mu_{t+s} \notin D_b \,\big\} \tag{3.46}$$

under the measure $\mathsf{P}_{t,x}$. Replacing $s$ with $\sigma_b$ in (3.36) and (3.41) and using the optional sampling theorem we find

$$\mathsf{E}_{t,x}\Big[V(t{+}\sigma_b, B^\mu_{t+\sigma_b})\Big] = V(t,x) \tag{3.47}$$

$$\begin{aligned} \mathsf{E}_{t,x}\Big[U^c(t{+}\sigma_b, B^\mu_{t+\sigma_b})\Big] &= U^c(t,x) \\ &\quad - \mathsf{E}_{t,x}\Big[\int_0^{\sigma_b} H(t{+}s, B^\mu_{t+s})\, I\big(B^\mu_{t+s} \notin D_c\big)\; ds\Big]. \end{aligned} \tag{3.48}$$

Since $(t,x)$ belongs to both $D_b$ and $D_c$ it follows that $U^c(t,x) = V(t,x) = 0$, and the fact that $U^c(t+\sigma_b, B^\mu_{t+\sigma_b}) \geq V(t+\sigma_b, B^\mu_{t+\sigma_b}) = 0$ implies

$$\mathsf{E}_{t,x}\left[\int_0^{\sigma_b} H(t+s, B^\mu_{t+s})\, I\big(B^\mu_{t+s} \notin D_c\big)\; ds\right] \leq 0\,. \tag{3.49}$$

The assumption that either $b_+(t) < c_+(t)$ or $b_-(t) > c_-(t)$ together with the continuity of the functions $c_-, c_+, b_-$ and $b_+$ means that there exists a small enough $t < u \leq T$ so that $b_+(s) < c_+(s)$ or $b_-(s) > c_-(s)$ for all $s \in [t,u]$. Consequently the $\mathsf{P}_{t,x}$ probability of $B^\mu$ spending a strictly positive amount of time (w.r.t. Lebesgue measure) in either of these regions is strictly positive. Combined with the fact that both $D_b$ and $D_c$ are contained in $\{H \geq 0\}$, this forces the expectation above to be strictly positive and provides a contradiction.

(iv) We show that $D_c = D_b$. Suppose that this is not the case so that $c_+(t) < b_+(t)$ or $c_-(t) > b_-(t)$ for some $t \in [0,T]$. Choose any point $x \in \big(c_+(t), b_+(t)\big)$ or $x \in \big(b_-(t), c_-(t)\big)$ and consider the stopping time

$$\tau_D = \inf\big\{\, s \in [0, T-t] \;\big|\; B^\mu_{t+s} \in D_b \,\big\} \tag{3.50}$$

under $\mathsf{P}_{t,x}$. Inserting $\tau_D$ in (3.36) and (3.41), taking $\mathsf{P}_{t,x}$ expectations and using the optional sampling theorem we obtain

$$\mathsf{E}_{t,x}\left[\int_0^{\tau_D} H(t+s, B^\mu_{t+s})\, ds\right] = V(t,x) \tag{3.51}$$

$$\mathsf{E}_{t,x}\Big[U^c\big(t+\tau_D, B^\mu_{t+\tau_D}\big)\Big] = U^c(t,x) - \mathsf{E}_{t,x}\left[\int_0^{\tau_D} H(t+s, B^\mu_{t+s})\, I\big(B^\mu_{t+s} \notin D_c\big)\; ds\right]. \tag{3.52}$$

Since $D_b$ is contained in $D_c$ and $U^c = 0$ on $D_c$ we must have $U^c\big(t+\tau_D, B^\mu_{t+\tau_D}\big) = 0$, and using the fact that $U^c \geq V$ we find that

$$\mathsf{E}_{t,x}\left[\int_0^{\tau_D} H(t+s, B^\mu_{t+s})\, I\big(B^\mu_{t+s} \in D_c\big)\; ds\right] \leq 0\,. \tag{3.53}$$

However as before the continuity of the boundaries $b_-, b_+, c_-$ and $c_+$ and the fact that $D_c \subseteq \{H \geq 0\}$ forces the expectation to be strictly positive and provides a contradiction. We therefore conclude that $c_+(t) = b_+(t)$ and $c_-(t) = b_-(t)$ for all $t \in [0,T]$. □

**Remark 1.** On closer inspection of the proof of Lemma 1 it is clear that a similar argument could be made for a broad class of Markov processes. Indeed, in order to reduce the optimal prediction problem to an equivalent optimal stopping problem, the only difficulty lies in computing the conditional distribution of $g$. For any *continuous* time-homogenous Markov process $X$ a simple application of the Markov property upon noting that $\max_{t \leq s \leq T} X_s = \big(\max_{0 \leq s \leq T-t} X_s\big) \circ \theta_t$ where $\theta_t$ is the shift operator shows that

$$\mathsf{P}_x(g < t \,|\, \mathcal{F}^X_t) = \begin{cases} \mathsf{P}_{X_t}\big(\min\limits_{0 \leq s \leq T-t} X_s > 0\big) & \text{if } \; X_t > 0 \\ \qquad\quad 0 & \text{if } \; X_t = 0 \\ \mathsf{P}_{X_t}\big(\max\limits_{0 \leq s \leq T-t} X_s < 0\big) & \text{if } \; X_t < 0 \end{cases} \tag{3.54}$$

for $t \in (0,T]$, and the problem reduces to computing the right-hand side of (3.54). Note that $g \le t$ in (2.7) and (2.8) can be replaced by $g < t$ so that the latter excludes the case $X_t = 0$ for $t \in [0,T]$ as used in (3.54).

For *discontinuous* Markov processes the problem is more involved as the representation (3.54) above is no longer valid. In this case one can formulate an analogous optimal prediction problem as follows. Let $X = (X_t)_{0 \le t \le T}$ be a standard Markov process defined on a filtered probability space $(\Omega, \mathcal{F}, (\mathcal{F}_t)_{0 \le t \le T}, \mathsf{P}_x)$ where $X_0 = x$ under $\mathsf{P}_x$ for $x \in I\!R$. Assume that 0 is regular for itself in the sense that $\mathsf{P}_0(\tau_0 = 0) = 1$ where $\tau_0 = \inf\{\, t \in (0,T] \mid X_t = 0 \,\}$, let $\ell = (\ell_t)_{0 \le t \le T}$ denote the local time of $X$ at 0 (see [1, p. 216]), and set

$$g = \inf\{\, t \in [0,T] \mid \ell_t = \ell_T \,\} \,. \tag{3.55}$$

Recalling that $t \mapsto \ell_t$ is continuous we see that $\ell_g = \ell_T$, and if $X$ is continuous then $g$ coincides with the last zero of $X$ before $T$ since $t \mapsto \ell_t$ increases only when $X_t = 0$ for $t \in [0,T]$. Inserting $g$ from (3.55) into (2.2) we obtain an analogous optimal prediction problem which can be reduced to a standard optimal stopping problem. Indeed, since $\ell$ is an additive functional it follows by the Markov property of $X$ that

$$\begin{aligned} \mathsf{P}_x(g \le t \,|\, \mathcal{F}_t) &= \mathsf{P}_x(\ell_t = \ell_T \,|\, \mathcal{F}_t) = \mathsf{P}_x(\ell_t = \ell_t + \ell_{T-t} \circ \theta_t \,|\, \mathcal{F}_t) \\ &= \mathsf{P}_x(\ell_{T-t} \circ \theta_t = 0 \,|\, \mathcal{F}_t) = \mathsf{P}_{X_t}(\ell_{T-t} = 0) = F(T-t, X_t) \end{aligned} \tag{3.56}$$

where we set $F(t,x) = \mathsf{P}_x(\ell_t = 0)$ for $(t,x) \in [0,T] \times I\!R$. Having (3.56) we can then proceed as in Lemma 1 above and define the value function $V$ as in (2.12) where $H(t,x) = 2F(T-t,x) - 1$ for $(t,x) \in [0,T] \times I\!R$. Note that $t \mapsto F(t,x)$ is decreasing so that $t \mapsto V(t,x)$ is decreasing and consequently if $(t,x)$ belongs to the stopping set $D$ then all points $(t+s,x)$ belong to $D$ for $s \in [0,T-t]$ with $x \in I\!R$ given and fixed. Note also that if $X$ is continuous then

$$F(T-t,x) = \mathsf{P}_x(\ell_{T-t} = 0) = \begin{cases} \mathsf{P}_x\big(\min\limits_{0 \le s \le T-t} X_s > 0\big) & \text{if } x > 0 \\ 0 & \text{if } x = 0 \\ \mathsf{P}_x\big(\max\limits_{0 \le s \le T-t} X_s < 0\big) & \text{if } x < 0 \end{cases} \tag{3.57}$$

for $t \in (0,T]$ and $x \in I\!R$, which is in agreement with (3.54) above. Finally, it is also evident that the previous considerations are by no means restricted to dimension one.

**Remark 2.** Note that there are several options if one wishes to compute numerical values for the value function $V$ and the optimal stopping boundaries $b_-$ and $b_+$. The first is through the equations (3.5) and (3.6)+(3.7), which will involve the solution of a coupled system of highly non–linear integral equations (this is the method applied to obtain Figure 1 above). Another is through the free-boundary problem (3.22)–(3.24) using a finite difference approach. One could also use a discrete-time discrete-space Markov chain approximation to $B^\mu$ and solve the analogous discrete problem using backward induction. The various methods each have their own peculiarities and strengths, and before trying to extract numbers from the problem, it is worth considering them all in order to choose the most appropriate one.

Jacques du Toit
School of Mathematics
The University of Manchester
Oxford Road
Manchester M13 9PL
United Kingdom
Jacques.Du-Toit@postgrad.manchester.ac.uk

Goran Peskir
School of Mathematics
The University of Manchester
Oxford Road
Manchester M13 9PL
United Kingdom
goran@maths.man.ac.uk

Albert Shiryaev
Steklov Mathematical Institute
Gubkina str. 8
119991 Moscow
Russia
albertsh@mi.ras.ru